# On the space of Type-2 interval with limit, continuity and differentiability of Type-2 interval-valued functions

**Md Sadikur Rahman, Ali Akbar Shaikh and Asoke Kumar Bhunia**

Department of mathematics, The University of Burdwan, Burdwan, India-713104.

mdsadikur.95@gmail.com, aliashaikh@math.buruniv.ac.in, akbhunia@math.buruniv.ac.in

**Abstract**
This paper deals with the new concept of interval whose both the bounds themselves are also intervals. We name this new type of interval as Type-2 interval. Here we have introduced Type-2 interval-valued function and its properties. To serve this purpose, we have defined a distance on the set of all Type-2 intervals, named as extended Moore distance for Type-2 intervals which is a metric on the set of all Type-2 intervals. Then we have shown that the space of Type-2 interval is a complete metric space with respect to extended Moore distance. Then we have introduced the concept of limit-continuity for Type-2 interval-valued function of single variable and also, we have derived some elementary properties of this concept. Subsequently, we have presented the idea of generalised Hukuhara difference on the set of Type-2 intervals. Finally, using this difference, we have defined gH-differentiability of Type-2 interval-valued function and discussed some of its properties.

**Keywords** Type-2 interval valued function; Type-2 interval Module, Type-2 interval sequence; Limit-Continuity of Type-2 interval valued function; Generalised Hukuhara difference, gH-differentiability

## 1. Introduction

Over the last few decades, researchers have considered various approaches, like stochastic approach, fuzzy approach, fuzzy-stochastic approach and interval approach to handle the uncertainty of a real-life problem. Among these approaches, interval approach is more significant. In interval approach, an imprecise parameter is represented by an interval. Already, the classical interval analysis has been introduced by some famous researchers like Moore [1], Hansen and Walster [2], to handle interval uncertainty. Also, we may refer to some excellent books on interval analysis written by Moore [1], Moore et al. [3] and Alefeld and Herzberger [4]. In an interval, the upper and lower bounds are fixed. But in reality, there are several situations in which an imprecise parameter can be represented by an interval where the upper and lower bounds may not be fixed. These bounds are flexible. These bounds may be expressed by intervals. So an interval can be expressed with the help of two intervals, one for lower bound and another for upper bound. This new type of interval is called Type-2 interval. In this paper, we have introduced some basic concepts of Type-2 interval mathematics, metric on $I_2(R)$, Type-2 interval valued function and its limit-continuity and differentiability to handle this type of uncertainty.

Here, first of all, we have defined the basic concepts of Type-2 intervals. Then we have discussed the Type-2 interval mathematics and Type-2 interval valued function of real variables. After that, we have defined extended Moore distance, Module [1] on $I_2(R)$; sequence in $I_2(R)$ and some of its properties. Next, we have shown that $I_2(R)$ is complete metric space with respect to extended Moore distance. Also, we have defined limit, continuity of a Type-2 interval valued function at a point and some of its elementary properties. After that we have developed some characterisation theorems on limit and continuity. Subsequently, we have derived generalized Hukuhara difference [5] for two Type-2 intervals. Then we have defined gH-differentiability along with some of its properties. Finally, we have drawn a fruitful conclusion from all the above-mentioned propositions.

## 2. Notation
(i) $N$: The set of all natural numbers.
(ii) $R$: The set of all real numbers.
(iii) $R^n$: The set of all ordered n-tuple of real numbers.
(iv) $I_1(R)$: The set of all Type-1 intervals of real numbers.
(v) $I_2(R)$: The set of all Type-2 intervals of real numbers of the form $[(\underline{a}_L, \bar{a}_L), (\underline{a}_U, \bar{a}_U)]$.

## 3. Some Basic concepts
The Type-2 interval denoted by $[(\underline{a}_L, \bar{a}_L), (\underline{a}_U, \bar{a}_U)]$ is defined in terms of Type-1 intervals given by
$[(\underline{a}_L, \bar{a}_L), (\underline{a}_U, \bar{a}_U)] = [a_L, a_U]$, where $a_L \in [\underline{a}_L, \bar{a}_L]$ and $a_U \in [\underline{a}_U, \bar{a}_U]$.

**Definition** Let $A = [(\underline{a}_L, \bar{a}_L), (\underline{a}_U, \bar{a}_U)]$, $B = [(\underline{b}_L, \bar{b}_L), (\underline{b}_U, \bar{b}_U)]$ be two Type-2 intervals. Then we say that $A$ is contained in $B$ which is denoted by $A \subseteq_2 B$ and is defined by
$A \subseteq_2 B \Leftrightarrow [\underline{a}_L, \bar{a}_U] \subseteq [\underline{b}_L, \bar{b}_U]$ and $[\bar{a}_L, \underline{a}_U] \subseteq [\bar{b}_L, \underline{b}_U]$
i.e., $A \subseteq_2 B \Leftrightarrow \underline{b}_L \leq \underline{a}_L, \bar{b}_L \leq \bar{a}_L, \underline{a}_U \leq \underline{b}_U$ and $\bar{a}_U \leq \bar{b}_U$.
$A \subset_2 B$ if $A \subseteq_2 B$ and $A \neq B$.

**Equality of Type-2 intervals**
Let $A = [(\underline{a}_L, \bar{a}_L), (\underline{a}_U, \bar{a}_U)]$, $B = [(\underline{b}_L, \bar{b}_L), (\underline{b}_U, \bar{b}_U)]$ be two Type-2 intervals.
Then $A = B$ if $A \subseteq_2 B$ and $B \subseteq_2 A$.

i.e., $A = B$ if $\underline{a}_L = \underline{b}_L, \bar{a}_L = \bar{b}_L, \underline{a}_U = \underline{b}_U, \bar{a}_U = \bar{b}_U$.

**1-Degenerate Type-2 interval**

Let $A = [(\underline{a}_L, \bar{a}_L),(\underline{a}_U, \bar{a}_U)]$ be a Type-2 interval. We say that $A$ is 1-degenerate if $\underline{a}_L = \bar{a}_L = \underline{a}_U = \bar{a}_U$. In this case, a Type-2 interval contains a single real number $a$. By convention we can identify a 1-degenerate Type-2 interval $[(a,a),(a,a)]$ with the real number $a$.

**2-Degenerate Type-2 interval**

Let $A = [(\underline{a}_L, \bar{a}_L),(\underline{a}_U, \bar{a}_U)]$ be a Type-2 interval. We say that $A$ is 2-degenerate if $\underline{a}_L = \bar{a}_L$, $\underline{a}_U = \bar{a}_U$. This 2-degenerate Type-2 interval can be identified with Type-1 interval $[a_L, a_U]$.

**Different operations of Type-2 interval**

Let $A = [(\underline{a}_L, \bar{a}_L),(\underline{a}_U, \bar{a}_U)]$ and $B = [(\underline{b}_L, \bar{b}_L),(\underline{b}_L, \bar{b}_L)]$ be two Type-2 intervals

**3.1 Addition** The addition of two Type-2 intervals A and B is defined as follows:
$$A + B = [(\underline{a}_L, \bar{a}_L),(\underline{a}_U, \bar{a}_U)] + [(\underline{b}_L, \bar{b}_L),(\underline{b}_U, \bar{b}_U)] = [(\underline{a}_L + \underline{b}_L, \bar{a}_L + \bar{b}_L),(\underline{a}_U + \underline{b}_U, \bar{a}_U + \bar{b}_U)]$$

**3.2 Subtraction** The subtraction between two Type-2 intervals A and B be defined as follows
$$A - B = [(\underline{a}_L, \bar{a}_L),(\underline{a}_U, \bar{a}_U)] - [(\underline{b}_L, \bar{b}_L),(\underline{b}_U, \bar{b}_U)] = [(\underline{a}_L - \bar{b}_U, \bar{a}_L - \underline{b}_U),(\underline{a}_U - \bar{b}_L, \bar{a}_U - \underline{b}_L)]$$

**Example 3.1** Let $A = [(-5,-2),(-1,3)], B = [(-3,1),(3,6)]$.
$A + B = [(-8,-1),(2,9)]$. $A - B = [(-11,-5),(-2,6)]$.

**3.3 Scalar multiplication**
$$\lambda.A = \lambda.[(\underline{a}_L, \bar{a}_L),(\underline{a}_U, \bar{a}_U)] = \begin{cases} [(\lambda\underline{a}_L, \lambda\bar{a}_L),(\lambda\underline{a}_U, \lambda\bar{a}_U)] \text{ if } \lambda \geq 0, \\ [(\lambda\bar{a}_U, \lambda\underline{a}_U),(\lambda\bar{a}_L, \lambda\underline{a}_L)] \text{ if } \lambda < 0. \end{cases}$$

**3.4 Multiplication**
$$AB = [(\min C, \min D),(\max D, \max C)]$$
where $C = \{\underline{a}_L \underline{b}_L, \underline{a}_L \bar{b}_U, \bar{a}_U \bar{b}_U, \bar{a}_U \underline{b}_L\}, D = \{\bar{a}_L \bar{b}_L, \underline{a}_U \bar{b}_L, \underline{a}_U \underline{b}_U, \underline{a}_U \bar{b}_L\}$

**Example 3.2** Let $A = [(-4,-1),(2,5)], B = [(-6,-3),(-1,3)]$.
  (i) For, $\lambda = 2$, $2A = [(-8,-2),(4,10)]$.
Again for $\lambda = -2$, $(-2)A = [(-10,-4),(2,8)]$.

(ii) $AB = \left[\left(\min C, \min D\right), \left(\max D, \max C\right)\right]$

*where* $C = \{24, -12, 15, -30\}, D = \{3, -6, -2, -6\}$

$= \left[(-30, -6), (3, 24)\right].$

## 3.5 Division

$\dfrac{A}{B} = A.\left(\dfrac{1}{B}\right)$ provided $0 \notin B$.

**Example 3.3** Let $A = \left[(-2, -1), (1, 3)\right], B = \left[(1, 2), (3, 4)\right]$.

Here $0 \notin B$, Hence $\dfrac{A}{B}$ is defined.

$\therefore \dfrac{A}{B} = \left[\left(-2, -\dfrac{1}{3}\right), \left(\dfrac{1}{2}, 3\right)\right].$

## 4. Type-2 interval valued function of real variable

A Type-2 interval valued function of real variable is a function $F_2 : D \to I_2(R)$, where $D(\neq \phi) \subseteq R$. It can be written as $F_2(x) = \left[\left(\underline{f}_L(x), \overline{f}_L(x)\right), \left(\underline{f}_U(x), \overline{f}_U(x)\right)\right]$, where $\underline{f}_L, \overline{f}_L, \underline{f}_U, \overline{f}_U : D \to R$ and $\underline{f}_L(x) \leq \overline{f}_L(x) \leq \underline{f}_U(x) \leq \overline{f}_U(x), \forall x \in D.$

**Example 4.1** Let $F_2(x) = \left[(x-1, x), (x+1, x+2)\right]$. Clearly $F_2(x)$ is a Type-2 interval valued function as $F_2 : R \to I_2(R)$

**Example 4.2** Let $f : R \to R$ be a real valued function and $C = \left[\left(\underline{c}_L, \overline{c}_L\right), \left(\underline{c}_U, \overline{c}_U\right)\right] \in I_2(R)$. Then $F_2(x) = Cf(x), x \in D$ is a Type-2 interval valued function as $F_2 : R \to I_2(R)$.

## 5. Sequence in $I_2(R)$

In this section, first of all we have defined the extended Moore distance on $I_2(R)$. Then we have defined module for Type-2 interval and introduced some properties of it. Also, we have defined sequence, convergent sequence, Cauchy sequence and presented some elementary properties of sequences in $I_2(R)$. Finally, we have shown that $I_2(R)$ is a complete metric space with respect to the Moore distance.

**Moore distance on $I_2(R)$**

Moore [1] introduced a distance on $I_1(R)$ which is named as Moore distance. Now we have generalised this distance on $I_2(R)$ and shown that it is a metric on $I_2(R)$. We have called this distance as extended Moore distance.

**Definition 5.1**

Let $X = [(\underline{x}_L, \bar{x}_L), (\underline{x}_U, \bar{x}_U)]$, $Y = [(\underline{y}_L, \bar{y}_L), (\underline{y}_U, \bar{y}_U)] \in I_2(R)$

Let us define the extended Moore distance between $X$ and $Y$ as

$$D_{M_2}(X,Y) = \max\{|\underline{x}_L - \underline{y}_L|, |\bar{x}_L - \bar{y}_L|, |\underline{x}_U - \underline{y}_U|, |\bar{x}_U - \bar{y}_U|\}.$$

**Theorem 5.1** The extended Moore distance defined earlier is a metric on $I_2(R)$.

**Proof:** (i) Let $X = [(\underline{x}_L, \bar{x}_L), (\underline{x}_U, \bar{x}_U)]$, $Y = [(\underline{y}_L, \bar{y}_L), (\underline{y}_U, \bar{y}_U)] \in I_2(R)$

From the definition, $D_{M_2}(X,Y) \geq 0$, and if $X = Y$, then $\underline{x}_L = \underline{y}_L, \bar{x}_L = \bar{y}_L, \underline{x}_U = \underline{y}_U, \bar{x}_U = \bar{y}_U$., therefore $D_{M_2}(X,Y) = 0$.

Conversely, if $D_{M_2}(X,Y) = 0$, then $|\underline{x}_L - \underline{y}_L| = |\bar{x}_L - \bar{y}_L| = |\underline{x}_U - \underline{y}_U| = |\bar{x}_U - \bar{y}_U| = 0$. Therefore, $\underline{x}_L = \underline{y}_L, \bar{x}_L = \bar{y}_L, \underline{x}_U = \underline{y}_U, \bar{x}_U = \bar{y}_U$ i.e., $X = Y$. (ii)

$$D_{M_2}(X,Y) = \max\{|\underline{x}_L - \underline{y}_L|, |\bar{x}_L - \bar{y}_L|, |\underline{x}_U - \underline{y}_U|, |\bar{x}_U - \bar{y}_U|\}$$
$$= \max\{|\underline{y}_L - \underline{x}_L|, |\bar{y}_L - \bar{x}_L|, |\underline{y}_U - \underline{x}_U|, |\bar{y}_U - \bar{x}_U|\} = D_{M_2}(Y,X).$$

(iii) Let $X = [(\underline{x}_L, \bar{x}_L), (\underline{x}_U, \bar{x}_U)]$, $Y = [(\underline{y}_L, \bar{y}_L), (\underline{y}_U, \bar{y}_U)]$, $Z = [(\underline{z}_L, \bar{z}_L), (\underline{z}_U, \bar{z}_U)] \in I_2(R)$

Now, $|\underline{x}_L - \underline{z}_L| \leq |\underline{x}_L - \underline{y}_L| + |\underline{y}_L + \underline{z}_L|$

$\leq \max\{|\underline{x}_L - \underline{y}_L|, |\bar{x}_L - \bar{y}_L|, |\underline{x}_U - \underline{y}_U|, |\bar{x}_U - \bar{y}_U|\} + \max\{|\underline{y}_L - \underline{z}_L|, |\bar{y}_L - \bar{z}_L|, |\underline{y}_U - \underline{z}_U|, |\bar{y}_U - \bar{z}_U|\}$

$= D_{M_2}(X,Y) + D_{M_2}(Y,Z)$.

Similarly, $|\bar{x}_L - \bar{z}_L| \leq D_{M_2}(X,Y) + D_{M_2}(Y,Z), |\underline{x}_U - \underline{z}_U| \leq D_{M_2}(X,Y) + D_{M_2}(Y,Z)$,

And $|\bar{x}_U - \bar{z}_U| \leq D_{M_2}(X,Y) + D_{M_2}(Y,Z)$

Thus, $D_{M_2}(X,Z) \leq D_{M_2}(X,Y) + D_{M_2}(Y,Z)$.

Therefore, $D_{M_2}$ is a metric on $I_2(R)$.

**Module for Type-2 intervals**

Here we have generalised the Module for Type-2 intervals as follows:

**Definition 5.2** Let $X \in I_2(R)$ be a Type-2 interval. The module of $X$ is defined as

$$\|X\|_{M_2} = D_{M_2}(X,0) = \max\{|\underline{x}_L|, |\bar{x}_L|, |\underline{x}_U|, |\bar{x}_U|\}.$$

**Module properties of Type-2 intervals**

The Module properties for Type-2 interval are as follows:

**Theorem 5.2**

(1). $\|X\|_{M_2} = 0 \Leftrightarrow X = 0$;

(2). $\|X + Y\|_{M_2} \leq \|X\|_{M_2} + \|Y\|_{M_2}$;

(3). $\|XY\|_{M_2} = \|X\|_{M_2} \|X\|_{M_2}$;

**Proof:** Proof follows from the definition

**Definition 5.3** The sequence of Type-2 intervals is a function $X: \mathbb{N} \to I_2(R)$. It is denoted by $\{X_n : n \in N\}$ or simply $\{X_n\}$.

**Example 5.1**

Let $X_n = \left[\left(\dfrac{1}{n+1}, \dfrac{1}{n}\right), \left(1+\dfrac{1}{n}, 2+\dfrac{1}{n}\right)\right] : n \in N$, Then $\{X_n\}$ be a sequence of Type-2 intervals.

**Definition 5.4** Let $\{X_n\}$ be a sequence in $I_2(R)$. Then $X_n$ is said to convergent to $X \in I_2(R)$ if for any $\varepsilon > 0, \exists\, n^o \in N$ st. $D_{M_2}(X_n, X) < \varepsilon, \forall n \geq n^o$, and it is denoted by $X_n \to X$ as $n \to \infty$.

**Theorem 5.3** Let $\left\{X_n = \left[\left(\underline{x}_L^{(n)}, \overline{x}_L^{(n)}\right), \left(\underline{x}_U^{(n)}, \overline{x}_U^{(n)}\right)\right]\right\}$ be a sequence in $I_2(R)$. Then $X_n \to X = \left[(\underline{x}_L, \overline{x}_L), (\underline{x}_U, \overline{x}_U)\right]$ if and only if $\underline{x}_L^{(n)} \to \underline{x}_L, \overline{x}_L^{(n)} \to \overline{x}_L, \underline{x}_U^{(n)} \to \underline{x}_U, \overline{x}_U^{(n)} \to \overline{x}_U$.

**Proof.** First let $X_n \to X = \left[(\underline{x}_L, \overline{x}_L), (\underline{x}_U, \overline{x}_U)\right]$.

Then for any $\varepsilon > 0, \exists N_o \in N$ st. $D_{M_2}(X_n, X) < \varepsilon, \forall n \geq N_o$.

i.e., $\max\left\{\left|\underline{x}_L^{(n)} - \underline{x}_L\right|, \left|\overline{x}_L^{(n)} - \overline{x}_L\right|, \left|\underline{x}_U^{(n)} - \underline{x}_U\right|, \left|\overline{x}_U^{(n)} - \overline{x}_U\right|\right\} < \varepsilon, \forall n \geq N_o$.

i.e., $\left|\underline{x}_L^{(n)} - \underline{x}_L\right| < \varepsilon, \left|\overline{x}_L^{(n)} - \overline{x}_L\right| < \varepsilon, \left|\underline{x}_U^{(n)} - \underline{x}_U\right| < \varepsilon, \left|\overline{x}_U^{(n)} - \overline{x}_U\right| < \varepsilon, \forall n \geq N_o$.

Hence $\underline{x}_L^{(n)} \to \underline{x}_L, \overline{x}_L^{(n)} \to \overline{x}_L, \underline{x}_U^{(n)} \to \underline{x}_U$ and $\overline{x}_U^{(n)} \to \overline{x}_U$, as $n \to \infty$.

Conversely let $\underline{x}_L^{(n)} \to \underline{x}_L, \overline{x}_L^{(n)} \to \overline{x}_L, \underline{x}_U^{(n)} \to \underline{x}_U, \overline{x}_U^{(n)} \to \overline{x}_U$ as $n \to \infty$.

Then for any $\varepsilon > 0$, $\exists\, N_1, N_2, N_3, N_4 \in N$ such that

$\left|\underline{x}_L^{(n)} - \underline{x}_L\right| < \varepsilon, \forall n \geq N_1$ \hfill (1)

$\left|\overline{x}_L^{(n)} - \overline{x}_L\right| < \varepsilon, \forall n \geq N_2$ \hfill (2)

$\left|\underline{x}_U^{(n)} - \underline{x}_U\right| < \varepsilon, \forall n \geq N_3$ \hfill (3)

$\left|\overline{x}_U^{(n)} - \overline{x}_U\right| < \varepsilon, \forall n \geq N_4$ \hfill (4)

Take, $N_o = \max\{N_1, N_2, N_3, N_4\}$.

Then from (1)-(4) we get,

$\left|\underline{x}_L^{(n)} - \underline{x}_L\right| < \varepsilon, \left|\overline{x}_L^{(n)} - \overline{x}_L\right| < \varepsilon, \left|\underline{x}_U^{(n)} - \underline{x}_U\right| < \varepsilon, \left|\overline{x}_U^{(n)} - \overline{x}_U\right| < \varepsilon, \forall n \geq N_o$.

i.e., $\max\left\{\left|\underline{x}_L^{(n)} - \underline{x}_L\right|, \left|\overline{x}_L^{(n)} - \overline{x}_L\right|, \left|\underline{x}_U^{(n)} - \underline{x}_U\right|, \left|\overline{x}_U^{(n)} - \overline{x}_U\right|\right\} < \varepsilon, \forall n \geq N_o$.

Since $\underline{x}_L^{(n)} \leq \overline{x}_L^{(n)} \leq \underline{x}_U^{(n)} \leq \overline{x}_U^{(n)}$, Then $\underline{x}_L \leq \overline{x}_L \leq \underline{x}_U \leq \overline{x}_U$.

Hence $D_{M_2}(X_n, X) < \varepsilon, \forall n \geq N_o$.

where $X_n = \left[\left(\underline{x}_L^{(n)}, \overline{x}_L^{(n)}\right), \left(\underline{x}_U^{(n)}, \overline{x}_U^{(n)}\right)\right], X = \left[(\underline{x}_L, \overline{x}_L), (\underline{x}_U, \overline{x}_U)\right]$.

Hence $X_n \to X$ as $n \to \infty$.

**Theorem 5.4** Let $\{X_n\}$ be a sequence in $I_2(R)$. If $X_n$ has a limit, it is unique.

**Proof.** In a metric space limit of a sequence if exists, it is unique. Since $(I_2(R), D_{M_2})$ be a metric space, $\{X_n\}$ has unique limit.

**Elementary Properties**

**Theorem 5.5** Let $\{X_n\}, \{Y_n\}$ be two sequences in $I_2(R)$ such that $X_n \to X, Y_n \to Y$, as $n \to \infty$. Then

(i) $X_n \pm Y_n \to X \pm Y$.

(ii) $X_n.Y_n \to X.Y$

(iii) $C.X_n \to C.X$, $C \in I_2(R)$.

(iv) $\dfrac{1}{Y_n} \to \dfrac{1}{Y}$, provided $0 \notin Y$.

(v) $\dfrac{X_n}{Y_n} \to \dfrac{X}{Y}$, provided $0 \notin X$.

**Proof.**
The proof of these properties are follows from the Theorem 5.2 and the elementary properties of real sequences.

**Definition** Let $\{X_n\}$ be a sequence in $I_2(R)$. Then it is said to be a Cauchy sequence if for any $\varepsilon > 0$, $\exists N_o \in N$ st. $D_{M_2}(X_n, X_m) < \varepsilon$, $\forall n, m \geq N_o$.

**Theorem 5.6** $(I_2, D_{M_2})$ is a complete metric space.

**Proof:** Let $X_n = \left\{\left[\left(\underline{x}_L^{(n)}, \overline{x}_L^{(n)}\right), \left(\underline{x}_U^{(n)}, \overline{x}_U^{(n)}\right)\right]\right\}$ be a Cauchy sequence in $I_2$.

Then for any $\varepsilon > 0$ $\exists n_o \in N$ such that $D_{M_2}(X_n, X_m) < \varepsilon, \forall n, m \geq n_o$.

$\Rightarrow \max\left\{\left|\underline{x}_L^{(n)} - \underline{x}_L^{(m)}\right|, \left|\overline{x}_L^{(n)} - \overline{x}_L^{(m)}\right|, \left|\underline{x}_U^{(n)} - \underline{x}_U^{(m)}\right|, \left|\overline{x}_U^{(n)} - \overline{x}_U^{(m)}\right|\right\} < \varepsilon, \forall n, m \geq n_o.$

$\Rightarrow \left|\underline{x}_L^{(n)} - \underline{x}_L^{(m)}\right| < \varepsilon, \left|\overline{x}_L^{(n)} - \overline{x}_L^{(m)}\right| < \varepsilon, \left|\underline{x}_U^{(n)} - \underline{x}_U^{(m)}\right| < \varepsilon, \left|\overline{x}_U^{(n)} - \overline{x}_U^{(m)}\right| < \varepsilon, \forall n, m \geq n_o$ \hfill (A)

$\Rightarrow \{\underline{x}_L^{(n)}\}, \{\overline{x}_L^{(n)}\}, \{\underline{x}_U^{(n)}\}$ and $\{\overline{x}_U^{(n)}\}$ all are Cauchy sequence in R.

Since R is complete,

$\exists$ unique $\underline{x}_L, \overline{x}_L, \underline{x}_U, \overline{x}_U$ in R such that $\underline{x}_L^{(n)} \to \underline{x}_L, \overline{x}_L^{(n)} \to \overline{x}_L, \underline{x}_U^{(n)} \to \underline{x}_U, \overline{x}_U^{(n)} \to \overline{x}_U$, as $n \to \infty$

Therefore in (A) taking $m \to \infty$, we get

$\left|\underline{x}_L^{(n)} - \underline{x}_L\right| < \varepsilon, \left|\overline{x}_L^{(n)} - \overline{x}_L\right| < \varepsilon, \left|\underline{x}_U^{(n)} - \underline{x}_U\right| < \varepsilon, \left|\overline{x}_U^{(n)} - \overline{x}_U\right| < \varepsilon, \forall n \geq n_o.$ \hfill (B)

Now since $x_n \leq y_n \Rightarrow \lim x_n \leq \lim y_n$. thus we get $\underline{x}_L \leq \overline{x}_L \leq \underline{x}_U \leq \overline{x}_U$. Therefore

$X = [(\underline{x}_L, \overline{x}_L), (\underline{x}_U, \overline{x}_U)] \in I_2$, Again from (B) we get $D_{M_2}(X_n, X) < \varepsilon, \forall n \geq n_o$.

Therefore $X_n \to X$ in $I_2$, as $n \to \infty$. Hence $(I_2, D_{M_2})$ is a complete metric space.

## 6.Limit, continuity and differentiability

In this section we have proposed the definition of limit and continuity of a Type-2 interval valued function of single variable at a point. Then we have developed some elementary properties, characterisation theorem for limit, continuity of a Type-2 interval valued function and limit, continuity for some special Type-2 interval valued function. Finally, we have presented the generalised Hukuhara difference [5] for Type-2 interval and gH-differentiability [5] of a Type-2 interval valued function and its characterisation theorem.

**Definition 6.1** Let $F_2 : D \subseteq R \to I_2(R)$ be a Type-2 interval valued function of single variable and $x_o$ be a limit point of $D$ and $L \in I_2(R)$. Then $F_2$ is said to tend to $L$ as $x$ tends to $x_o$, if for any $\varepsilon > 0$, $\exists \delta > 0$ such that $D_{M_2}(F_2(x), L) < \varepsilon$ whenever $0 < |x - x_o| < \delta$. And it is denoted by $\lim_{x \to x_o} F_2(x) = L$.

**Theorem 6.1** Let $F_2 : D \subseteq R \to I_2(R)$ be a Type-2 interval valued function such that $F(x) = \left[\left(\underline{f}_L(x), \overline{f}_L(x)\right), \left(\underline{f}_U(x), \overline{f}_U(x)\right)\right]$ and $x_o$ be a limit point of $D$. Then $\lim_{x \to x_o} F_2(x)$ exists iff $\lim_{x \to x_o} \underline{f}_L(x), \lim_{x \to x_o} \overline{f}_L(x), \lim_{x \to x_o} \underline{f}_U(x)$ and $\lim_{x \to x_o} \overline{f}_U(x)$ exist. And

$$\lim_{x \to x_o} F_2(x) = \left[\left(\lim_{x \to x_o} \underline{f}_L(x), \lim_{x \to x_o} \overline{f}_L(x)\right), \left(\lim_{x \to x_o} \underline{f}_U(x), \lim_{x \to x_o} \overline{f}_U(x)\right)\right]$$

**Proof.**

Let $\lim_{x \to x_o} F_2(x)$ exist and $\lim_{x \to x_o} F_2(x) = L = \left[\left(\underline{l}_L, \overline{l}_L\right), \left(\underline{l}_U, \overline{l}_U\right)\right]$.

Then for any $\varepsilon > 0, \exists \delta > 0$ such that $D_{M_2}(F_2(x), L) < \varepsilon$ whenever $0 < |x - x_o| < \delta$.

i.e., $\max\left\{\left|\underline{f}_L(x) - \underline{l}_L\right|, \left|\overline{f}_L(x) - \overline{l}_L\right|, \left|\underline{f}_U(x) - \underline{l}_U\right|, \left|\overline{f}_U(x) - \overline{l}_L\right|\right\} < \varepsilon$ whenever $0 < |x - x_o| < \delta$.

i.e., $\left|\underline{f}_L(x) - \underline{l}_L\right| < \varepsilon, \left|\overline{f}_L(x) - \overline{l}_L\right| < \varepsilon, \left|\underline{f}_U(x) - \underline{l}_U\right| < \varepsilon, \left|\overline{f}_U(x) - \overline{l}_U\right| < \varepsilon$ whenever $0 < |x - x_o| < \delta$.

i.e., $\lim_{x \to x_0} \underline{f}_L(x) = \underline{l}_L, \lim_{x \to x_o} \overline{f}_L(x) = \overline{l}_L, \lim_{x \to x_0} \underline{f}_U(x) = \underline{l}_U, \lim_{x \to x_o} \overline{f}_U(x) = \overline{l}_U$.

Hence all $\lim_{x \to x_o} \underline{f}_L(x), \lim_{x \to x_o} \overline{f}_L(x), \lim_{x \to x_o} \underline{f}_U(x), \lim_{x \to x_o} \overline{f}_U(x)$ exist.

And $\lim_{x \to x_o} F_2(x) = \left[\left(\underline{l}_L, \overline{l}_L\right), \left(\underline{l}_U, \overline{l}_U\right)\right] = \left[\left(\lim_{x \to x_o} \underline{f}_L(x), \lim_{x \to x_o} \overline{f}_L(x)\right), \left(\lim_{x \to x_o} \underline{f}_U(x), \lim_{x \to x_o} \overline{f}_U(x)\right)\right]$

Conversely let $\lim_{x \to x_o} \underline{f}_L(x), \lim_{x \to x_o} \overline{f}_L(x), \lim_{x \to x_o} \underline{f}_U(x)$ and $\lim_{x \to x_o} \overline{f}_U(x)$ exist.

Let $\lim_{x \to x_0} \underline{f}_L(x) = \underline{l}_L, \lim_{x \to x_o} \overline{f}_L(x) = \overline{l}_L, \lim_{x \to x_0} \underline{f}_U(x) = \underline{l}_U, \lim_{x \to x_o} \overline{f}_U(x) = \overline{l}_U$.

Then for any $\varepsilon > 0, \exists \delta_1, \delta_2, \delta_3, \delta_4 > 0$ such that

$$|\underline{f}_L - \underline{l}_L| < \varepsilon \text{ whenever } 0 < |x - x_o| < \delta_1, \tag{1}$$

$$|\overline{f}_L - \overline{l}_L| < \varepsilon \text{ whenever } 0 < |x - x_o| < \delta_2, \tag{2}$$

$$|\underline{f}_U - \underline{l}_U| < \varepsilon \text{ whenever } 0 < |x - x_o| < \delta_3, \tag{3}$$

$$|\overline{f}_U - \overline{l}_U| < \varepsilon \text{ whenever } 0 < |x - x_o| < \delta_4, \tag{4}$$

Let us consider $\delta = \min\{\delta_1, \delta_2, \delta_3, \delta_4\}$.

Then from (1)-(4) we get,

$|\underline{f}_L(x) - \underline{l}_L| < \varepsilon, |\overline{f}_L(x) - \overline{l}_L| < \varepsilon, |\underline{f}_U(x) - \underline{l}_U| < \varepsilon, |\overline{f}_U(x) - \overline{l}_U| < \varepsilon$ whenever $0 < |x - x_o| < \delta$.

i.e., $\max\{|\underline{f}_L(x) - \underline{l}_L|, |\overline{f}_L(x) - \overline{l}_L|, |\underline{f}_U(x) - \underline{l}_U|, |\overline{f}_U(x) - \overline{l}_L|\} < \varepsilon$ whenever $0 < |x - x_o| < \delta$. (5)

Since $\underline{f}_L(x) \leq \overline{f}_L(x) \leq \underline{f}_U(x) \leq \overline{f}_U(x)$, $\underline{l}_L \leq \overline{l}_L \leq \underline{l}_U \leq \overline{l}_U$.

Now let $L = \left[(\underline{l}_L, \overline{l}_L), (\underline{l}_U, \overline{l}_U)\right]$.

Then from (5) we get,

$D_{M_2}(F_2(x), L) < \varepsilon$ whenever $0 < |x - x_o| < \delta$.

Hence $\lim_{x \to x_o} F_2(x) = L = \left[(\underline{l}_L, \overline{l}_L), (\underline{l}_U, \overline{l}_U)\right]$.

**Theorem 6.2** Let $F_2 : D \subseteq R \to I_2(R)$ with $F_2(x) = Cf(x)$, where $C = \left[(\underline{c}_L, \overline{c}_L), (\underline{c}_U, \overline{c}_U)\right]$ and $f : D \to R$ be real valued function. Let $x_o$ be a limit point of $D$. If $\lim_{x \to x_o} f(x)$ exists, then $\lim_{x \to x_o} F_2(x)$ exists and $\lim_{x \to x_o} F_2(x) = C. \lim_{x \to x_o} f(x)$.

**Proof.** Assume that $\lim_{x \to x_o} f(x) = l$. and let $c_k = \|C\|_{M_2} > 0$.

If possible, let $\lim_{x \to x_o} F_2(x) \neq C. \lim_{x \to x_o} f(x)$. Then

$\exists \varepsilon > 0$, such that $\forall \delta > 0$ and $\forall x : 0 < |x - x_o| < \delta$, $D_{M_2}(C.f(x), Cl) \geq \varepsilon$.

So, $\max\{|\underline{c}_L f(x) - \underline{c}_L l|, |\overline{c}_L f(x) - \overline{c}_L l|, |\underline{c}_U f(x) - \underline{c}_U l|, |\overline{c}_U f(x) - \overline{c}_U l|\} \geq \varepsilon$.

i.e., $\max\{|\underline{c}_L|, |\overline{c}_L|, |\underline{c}_U|, |\overline{c}_U|\} \cdot |f(x) - l| \geq \varepsilon$.

i.e., $\|C\|_{M_2} \cdot |f(x) - l| \geq \varepsilon$,

i.e., $c_k \cdot |f(x) - l| \geq \varepsilon$, $\forall x : 0 < |x - x_o| < \delta$.

This implies, $|f(x) - l| \geq \dfrac{\varepsilon}{c_k}$, $\forall x : 0 < |x - x| < \delta$.

i.e., $\lim_{x \to x_o} f(x) = l$, which contradicts our assumption.

Therefore, $\lim_{x \to x_o} F_2(x) = C. \lim_{x \to x_o} f(x)$.

If $c_k = 0$, then clearly $\lim_{x \to x_o} F_2(x) = C. \lim_{x \to x_o} f(x)$.

**Note 6.1** Converse of the above theorem is not true.

For example, let us consider the function $f : R \to R$ defined by $f(x) = \begin{cases} 1, & \text{if } x \geq 0 \\ -1, & \text{if } x \leq 0 \end{cases}$

It is clear that $\lim_{x \to x_o} f(x)$ does not exist.

Let us define $F_2(x)$ by $F_2(x) = Cf(x)$, where $C = [(-2,-1),(1,2)]$.

Hence $F_2(x) = [(-2,-1),(1,2)], \forall x \in R.$

So $\lim_{x \to 0} F_2(x) = [(-2,-1),(1,2)]$.

**Theorem 6.3** Let $F_2 : D \to I_2(R)$ with $F_2(x) = Cf(x), C = [(\underline{c}_L, \overline{c}_L),(\underline{c}_U, \overline{c}_U)]$ and $x_o$ be the limit point of $D$. If $\lim_{x \to x_o} F_2(x)$ exists, then one of the following cases holds:

(a) There exists $\lim_{x \to x_o} f(x)$ and $\lim_{x \to x_o} F_2(x) = C \lim_{x \to x_o} f(x)$.

(b) There exists $\lim_{x \to x_o} |f(x)|; [\underline{c}_L, \overline{c}_L] = -[\underline{c}_U, \overline{c}_U]$.

**Proof.** Here $F_2$ can be rewritten as

$$F_2(x) = [(m_1(x), m_2(x)),(m_3(x), m_4(x))]$$

where

$m_1(x) = \min\{\underline{c}_L g(x), \overline{c}_U g(x)\}$

$m_2(x) = \min\{\overline{c}_L g(x), \underline{c}_U g(x)\}$

$m_3(x) = \max\{\overline{c}_L g(x), \underline{c}_U g(x)\}$

$m_4(x) = \max\{\underline{c}_L g(x), \overline{c}_U g(x)\}$

Since $\lim_{x \to x_0} F_2(x)$ exists, then by Theorem 2.1, the limits $\lim_{x \to x_o} m_i(x)$ exist, $i = 1, 2, 3, 4$.

So $\lim_{x \to x_o}(m_1(x) + m_4(x))$ and $\lim_{x \to x_o}(m_2(x) + m_3(x))$ also exist.

But $m_1(x) + m_4(x) = (\underline{c}_L + \overline{c}_U)g(t)$ and $m_2(x) + m_3(x) = (\overline{c}_L + \underline{c}_U)g(x)$.

If $(\underline{c}_L + \overline{c}_U) \neq 0$ or $(\overline{c}_L + \underline{c}_U) \neq 0$, then $\lim_{x \to x_o} g(x)$ exists. Hence (a) holds.

If $\underline{c}_L + \overline{c}_U = 0 = \overline{c}_L + \underline{c}_U$ then $\underline{c}_L = -\overline{c}_U, \overline{c}_L = -\underline{c}_U$.

and hence

$[\underline{c}_L, \overline{c}_L] = [-\overline{c}_U, -\underline{c}_U] = -[\underline{c}_U, \overline{c}_U]$ and

$F(x) = [(-\overline{c}_U, -\underline{c}_U),(\underline{c}_U, \overline{c}_U)] \cdot g(x) = [(-\overline{c}_U |g(x)|, -\underline{c}_U |g(x)|),(\underline{c}_U |g(x)|, \overline{c}_U |g(x)|)]$

Then from Theorem 2.1, we can say that $\lim_{x \to x_o} |g(x)|$ exists. Hence proves the theorem 2.3(b).

**Definition 6.2** Let $F_2 : D \subseteq R \to I_2(R)$ be a Type-2 interval valued function of single variable and $x_o \in D$. and $L \in I_2(R)$. Then $F_2$ is said to be continuous at $x = x_o$, if for any

$\varepsilon > 0, \exists \delta > 0$ such that. $D_{M_2}(F_2(x), F_2(x_o)) < \varepsilon$ whenever $|x - x_o| < \delta$. i.e., $F_2$ is continuous at $x = x_o$ if $\lim_{x \to x_o} F_2(x) = F_2(x_o)$.

**Theorem 6.4** Let $F_2 : D \to I_2(R)$ be given by $F_2(x) = \left[\left(\underline{f}_L(x), \overline{f}_L(x)\right), \left(\underline{f}_U(x), \overline{f}_U(x)\right)\right]$. Then $F_2$ is continuous at $x = x_o \in D$, iff $\underline{f}_L, \overline{f}_L, \underline{f}_U, \overline{f}_U$ all are continuous at $x = x_o$.

**Theorem 6.5** Let $F_2 : D \subseteq R \to I_2(R)$ with $F_2(x) = Cf(x), C = \left[(\underline{c}_L, \overline{c}_L), (\underline{c}_U, \overline{c}_U)\right]$, and $f : D \to R$ be real valued function. Let $x_o \in D$. If $f$ continuous at $x = x_o$, then $F_2$ is continuous at $x = x_o$.

**Theorem 6.6** Let $F_2 : D \to I_2(R)$ with $F_2(x) = Cf(x), C = \left[(\underline{c}_L, \overline{c}_L), (\underline{c}_U, \overline{c}_U)\right]$. Let $x_o \in D$. If $F_2$ is continuous at $x = x_o$, then one of the following cases holds:
(a) $f$ is continuous at $x = x_o$.
(b) $|f|$ is continuous at $x_o$ and $[\underline{c}_L, \overline{c}_L] = -[\underline{c}_U, \overline{c}_U]$.

**Generalised Hukuhara difference of Type-2 intervals:**

The generalised Hukuhara difference [5] for two non-empty compact subsets of $R^n$ have already been introduced. Here we have introduced the generalised Hukuhara difference for any nonempty sets and obtained the generalised Hukuhara difference for two non-empty sets in $I_2(R)$.

**Definition [5]** Let $A$ and $B$ be two non-empty sets. Then the generalised Hukuhara difference $A \ominus_g B$ of $A$ and $B$ is defined by $A \ominus_g B = C \Leftrightarrow \begin{cases} (a) A = B + C \\ \text{or} \\ (b) B = A + (-1)C \end{cases}$

**Proposition 6.1** Let $A = \left[(\underline{a}_L, \overline{a}_L), (\underline{a}_U, \overline{a}_U)\right], B = \left[(\underline{b}_L, \overline{b}_L), (\underline{b}_U, \overline{b}_U)\right] \in I_2(R)$. Then $A \ominus_g B$ is given by $A \ominus_g B = C = \left[(\underline{c}_L, \overline{c}_L), (\underline{c}_U, \overline{c}_U)\right]$,
where
$\underline{c}_L = \min\{\underline{a}_L - \underline{b}_L, \overline{a}_U - \overline{b}_U\}; \overline{c}_U = \max\{\underline{a}_L - \underline{b}_L, \overline{a}_U - \overline{b}_U\};$
$\overline{c}_L = \min\{\overline{a}_L - \overline{b}_L, \underline{a}_U - \underline{b}_U\}; \underline{c}_U = \max\{\overline{a}_L - \overline{b}_L, \underline{a}_U - \underline{b}_U\}.$

**Proof.** From the definition of generalised Hukuhara difference we get,

$A \ominus_g B = C \Leftrightarrow \begin{cases} (a) A = B + C \\ \text{or} \\ (b) B = A + (-1)C \end{cases}$

Now,

$A = B + C$

$\Leftrightarrow [(\underline{a}_L, \bar{a}_L), (\underline{a}_U, \bar{a}_U)] = [(\underline{b}_L + \underline{c}_L, \bar{b}_L + \bar{c}_L), (\underline{b}_U + \underline{c}_U, \bar{b}_U + \bar{c}_U)]$

$\Leftrightarrow \underline{a}_L = \underline{b}_L + \underline{c}_L, \bar{a}_L = \bar{b}_L + \bar{c}_L, \underline{a}_U = \underline{b}_U + \underline{c}_U, \bar{a}_U = \bar{b}_U + \bar{c}_U.$

$\Leftrightarrow \underline{c}_L = \underline{a}_L - \underline{b}_L, \bar{c}_L = \bar{a}_L - \bar{b}_L, \underline{c}_U = \underline{a}_U - \underline{b}_U, \bar{c}_U = \bar{a}_U - \bar{b}_U.$

Again,

$B = A + (-1)C$

$\Leftrightarrow [(\underline{b}_L, \bar{b}_L), (\underline{b}_U, \bar{b}_U)] = [(\underline{a}_L, \bar{a}_L), (\underline{a}_U, \bar{a}_U)] + (-1)[(\underline{c}_L, \bar{c}_L), (\underline{c}_U, \bar{c}_U)].$

$\Leftrightarrow [(\underline{b}_L, \bar{b}_L), (\underline{b}_U, \bar{b}_U)] = [(\underline{a}_L - \bar{c}_U, \bar{a}_L - \underline{c}_U), (\underline{a}_U - \bar{c}_L, \bar{a}_U - \underline{c}_L)].$

$\Leftrightarrow \underline{b}_L = \underline{a}_L - \bar{c}_U, \bar{b}_L = \bar{a}_L - \underline{c}_U, \underline{b}_U = \underline{a}_U - \bar{c}_L, \bar{b}_U = \bar{a}_U - \underline{c}_L.$

$\Leftrightarrow \bar{c}_U = \underline{a}_L - \underline{b}_L, \underline{c}_U = \bar{a}_L - \bar{b}_L, \bar{c}_L = \underline{a}_U - \underline{b}_U, \underline{c}_L = \bar{a}_U - \bar{b}_U.$

Hence,

$\underline{c}_L = \min\{\underline{a}_L - \underline{b}_L, \bar{a}_U - \bar{b}_U\}; \bar{c}_U = \max\{\underline{a}_L - \underline{b}_L, \bar{a}_U - \bar{b}_U\};$

$\bar{c}_L = \min\{\bar{a}_L - \bar{b}_L, \underline{a}_U - \underline{b}_U\}; \underline{c}_U = \max\{\bar{a}_L - \bar{b}_L, \underline{a}_U - \underline{b}_U\}.$

**Definition 6.3** Let $F_2 : D \subseteq R \to I_2(R)$ be a Type-2 interval valued function and $x_o \in D.$ Then $F_2$ is said to be gH-differentiable at $x_o$, if the following limit

$\lim_{h \to 0} \dfrac{F_2(x_o + h) - F_2(x_o)}{h}$ exists. where $h$ be such that $x_o + h \in D.$

It is denoted by $F_2'(x_o) = \lim_{h \to 0} \dfrac{F_2(x_o + h) - F_2(x_o)}{h}.$

**Theorem 6.7** Let $F_2 : D \to I_2(R)$ be a Type-2 interval valued function such that $F_2(x) = [(\underline{f}_L(x), \bar{f}_L(x)), (\underline{f}_U(x), \bar{f}_U(x))].$ Then $F_2$ is gH-differentiable at $x_o \in D$ iff the following condition holds:

$\underline{f}_L, \bar{f}_L, \underline{f}_U$ and $\bar{f}_U$ are differentiable at $x_o$ and

$F_2'(x_o) = [(\min C, \min D), (\max D, \max C)],$

where

$C = \{\underline{f}_L'(x_o), \bar{f}_U'(x_o)\};$

$D = \{\bar{f}_L'(x_o), \underline{f}_U'(x_o)\}.$

**Proof.**

Since $F_2(x)$ be gH-differentiable at $x_o \in D.$

Then

$\lim_{h \to 0} \dfrac{F_2(x_o + h) \ominus_{gH} F_2(x_o)}{h}$ exists

and $F_2'(x_o) = \lim_{h \to 0} \dfrac{F_2(x_o + h) \ominus_{gH} F_2(x_o)}{h} \in I_2(R).$

Now $F(x_o + h) \ominus_{gH} F(x_o) = [(\min P, \min Q), (\max Q, \max P)]$

where,

$P = \{(\underline{f_L}(x_o + h) - \underline{f_L}(x_o)), (\overline{f_U}(x_o + h) - \overline{f_U}(x_o))\};$

$Q = \{(\overline{f_L}(x_o + h) - \overline{f_L}(x_o)), (\underline{f_U}(x_o + h) - \underline{f_U}(x_o))\}.$

Since $\lim_{h \to 0} \dfrac{F_2(x_o + h) \ominus_{gH} F_2(x_o)}{h}$ exists iff

$\lim_{h \to 0} \dfrac{\min P}{h}, \lim_{h \to 0} \dfrac{\min Q}{h}, \lim_{h \to 0} \dfrac{\max Q}{h}$ and $\lim_{h \to 0} \dfrac{\max P}{h}$ exist.

i.e., iff

$\underline{f_L}'(x_o) = \lim_{h \to 0} \dfrac{\underline{f_L}(x_o + h) - \underline{f_L}(x_o)}{h}, \overline{f_L}'(x_o) = \lim_{h \to 0} \dfrac{\overline{f_L}(x_o + h) - \overline{f_L}(x_o)}{h},$

$\underline{f_U}'(x_o) = \lim_{h \to 0} \dfrac{\underline{f_U}(x_o + h) - \underline{f_U}(x_o)}{h}$ and $\overline{f_U}'(x_o) = \lim_{h \to 0} \dfrac{\overline{f_U}(x_o + h) - \overline{f_U}(x_o)}{h}$ exist.

i.e., iff $\underline{f_L}, \overline{f_L}, \underline{f_U}, \overline{f_U}$ all are differentiable at $x_o$.

From the Theorem 2.1, we get,

$$F_2'(x_o) = \begin{bmatrix} (\min\{\underline{f_L}'(x_o), \overline{f_U}'(x_o)\}, \min\{\overline{f_L}'(x_o), \underline{f_U}'(x_o)\}), \\ (\max\{\overline{f_L}'(x_o), \underline{f_U}'(x_o)\}, \max\{\underline{f_L}'(x_o), \overline{f_U}'(x_o)\}) \end{bmatrix}$$

**Definition 6.4** Let $F_2 : D \to I_2(R)$ be a Type-2 interval valued function such that $F_2(x) = [(\underline{f_L}(x), \overline{f_L}(x)), (\underline{f_U}(x), \overline{f_U}(x))]$. Then

(i) $F_2$ is said to be gH-differentiable at $x_o \in D$ in first form if $\underline{f_L}, \overline{f_L}, \underline{f_U}$ and $\overline{f_U}$ are differentiable at $x_o$ and $F_2'(x_o) = [(\underline{f_L}'(x_o), \overline{f_L}'(x_o)), (\underline{f_U}'(x_o), \overline{f_U}'(x_o))]$.

(ii) $F_2$ is said to be gH-differentiable at $x_o \in D$ in second form if $\underline{f_L}, \overline{f_L}, \underline{f_U}$ and $\overline{f_U}$ are differentiable at $x_o$ and $F_2'(x_o) = [(\overline{f_U}'(x_o), \underline{f_U}'(x_o)), (\overline{f_L}'(x_o), \underline{f_L}'(x_o))]$.

**Theorem 6.8** Let $F_2 : D \to I_2(R)$ be a bounded Type-2 interval valued function such that $F_2(x) = Cf(x)$ with $C = [(\underline{c_L}, \overline{c_L}), (\underline{c_U}, \overline{c_U})]$.

Then $F_2$ is gH-differentiable at $x_o \in D$, if $f$ is differentiable at $x_o \in D$ and has the same sign throughout D then

$F_2'(x_o) = Cf'(x_o).$

**Proof.** Let $f(x) \geq 0$,

$\therefore F_2(x) = Cf(x) = [(\underline{c_L} f(x), \overline{c_L} f(x)), (\underline{c_U} f(x), \overline{c_U} f(x))].$

Since $f(x)$ is differentiable at $x_o$, then $\underline{c}_L f(x), \overline{c}_L f(x), \underline{c}_U f(x), \overline{c}_U f(x)$ all are differentiable at $x_o$. Then by the Theorem 2.7 $F_2(x)$ is gH-differentiable and

$$F_2'(x_o) = \begin{bmatrix} \left(\min\{\underline{c}_L f'(x_o), \overline{c}_U f'(x_o)\}, \min\{\overline{c}_L f'(x_o), \underline{c}_U f'(x_o)\}\right), \\ \left(\max\{\overline{c}_L f'(x_o), \underline{c}_U f'(x_o)\}, \max\{\underline{c}_L f'(x_o), \overline{c}_U f'(x_o)\}\right) \end{bmatrix}$$
$$= \left[(\underline{c}_L, \overline{c}_L), (\underline{c}_U, \overline{c}_L)\right] \cdot \left[(f'(x_o), f'(x_o)), (f'(x_o), f'(x_o))\right]$$
$$= Cf'(x_o).$$

If $f(x) \leq 0$, the proof is similar as above.

## 7. Conclusion

In this paper, for first time, the analytical concepts of Type-2 intervals have been introduced. Then the definitions of metric, module, sequence of Type-2 interval numbers and, limit, continuity and gH-differentiability of Type-2 interval-valued functions have been proposed. For further works, based on the concepts developed in this paper, one may introduce uniform convergence and the integrability of Type-2 interval-valued function; initial-value problem of ordinary differential equation and its solution methodology. Also, the concept of Type-2 intervals and interval-valued functions may be extended to the theory of optimization and applications.